\documentclass[11pt]{amsart}

\usepackage{fullpage,amssymb}

\newcommand{\Z}{\mathbb{Z}}

\newcommand{\Q}{\mathbb{Q}}

\newcommand{\ZN}[1]{\Z/{#1}\Z}

\newcommand{\K}{\mathbf{K}}
\newcommand{\OK}{\mathcal{O}_K}
\newcommand{\KH}{\mathbf{K}_H}

\newcommand {\legendre}[2]{\genfrac {(}{)}{1pt}{}{#1}{#2}}

\newtheorem{proposition}{Proposition}[section]

\newif\ifdraft \drafttrue

\begin{document}

\ifdraft
\title[fastECPP]{Implementing the asymptotically fast version of the
elliptic curve primality proving algorithm}
%{\bf version 050127}}
\fi
\author{F. Morain}
\address[F. Morain]{LIX \'Ecole Polytechnique, CNRS/UMR 7161, INRIA/Futurs,
F-91128 Palaiseau CEDEX, FRANCE}
\email[F. Morain]{morain@lix.polytechnique.fr}
\thanks{Projet TANC, P\^ole Commun
de Recherche en Informatique du Plateau de Saclay, CNRS, \'Ecole
polytechnique, INRIA, Universit\'e Paris-Sud. The author is on leave
from the French Department of Defense, D\'el\'egation G\'en\'erale
pour l'Armement.}
\date{\today}

\begin{abstract}
The elliptic curve primality proving (ECPP) algorithm is one of the current
fastest practical algorithms for proving the primality of large
numbers. Its running time cannot be proven rigorously, but heuristic
arguments show that it should run in time $\tilde{O}((\log N)^5)$ to
prove the primality of $N$. An asymptotically fast version of it,
attributed to J.~O.~Shallit, runs in time $\tilde{O}((\log N)^4)$. The
aim of this article is to describe this version in more details,
leading to actual implementations able to handle numbers with
several thousands of decimal digits.
\end{abstract}

\maketitle

\ifdraft
\else
\tableofcontents
\fi

%\input{fastecpp-body}
%%%%% S
\section{Introduction}

From the work of Agrawal, Kayal and Saxena \cite{AgKaSa02}, we know
that determining
the primality of an integer $N$ can be done in proven deterministic
polynomial time $\tilde{O}((\log N)^{10.5})$. More recently, H.-W.~Lenstra and
C.~Pomerance have announced a version in $\tilde{O}((\log
N)^{6})$. Building on the work of P.~Berrizbeitia
\cite{Berrizbeitia02}, D.~Bernstein \cite{Bernstein03} and
P.~Mih\u{a}ilescu \& R.~Mocenigo \cite{MiAv03}, independently, have
given improved probabilistic versions with a claim of proven complexity of
$\tilde{O}((\log N)^{4})$, reusing classical cyclotomic ideas that
originated in the Jacobi sums test \cite{AdPoRu83,CoLe84}. For more on
primality before AKS, we refer the reader to \cite{CrPo00}. For the
recent developments, see \cite{Bernstein03a}.

All the known versions of the AKS algorithm are for the time being too
slow to prove the primality of large explicit numbers. 
On the other hand, the elliptic curve primality proving algorithm
\cite{AtMo93b} has been used for years to prove the primality of
always larger numbers\footnote{See the web page of M.~Martin, {\tt
http://www.ellipsa.net/}, or that of the author}.
The algorithm has a heuristic running time of
$\tilde{O}((\log N)^5)$. In the
course of writing \cite{Morain04}, the author rediscovered the article
\cite{LeLe90}, in which an asymptotically fast version of ECPP is
described. This version, attributed to J.~O.~Shallit, has a heuristic
running time of $\tilde{O}((\log N)^{4})$.
\iffalse
To the best of the author's knowledge, this fast
version has never been implemented. Due to the large numbers now being
dealt with and the stimulating effect of both the AKS thread and the
fair-play concurrence of PRIMO, the author realized that
implementing this fast version was now at hand, helped by recent work
on the computation of defining polynomials of Hilbert class fields. 
\fi
The aim of this paper is to describe {\sc fastECPP}, give
a heuristic analysis of it and describe its implementation. 

Section 2 collects some well-known facts on imaginary quadratic
fields, that can be found for instance in \cite{Cox89}. Section 3
presents the basic ECPP algorithm and analyzes
it. In Section 4, the fast version is described and its complexity
estimated. Section 5 explains the implementation and Section 6 gives
some actual timings on large numbers.
\iffalse
Section 7 describes a compact fingerprint that can be used to rebuild
a certificate from it.
\fi

%%%%% S
\section{Quadratic fields}

A discriminant $-D<0$ is said to be fundamental if and only if $D$ is free of
odd square prime factors, and moreover $D \equiv 3\bmod 4$ or when
$4\mid D$, $(D/4) \bmod 4 \in \{1, 2\}$. The quantity
$$\mathcal{D}(X) = \# \{D \leq X, -D \text{ is fundamental} \}$$
is easily seen to be $O(X)$.

A fundamental discriminant may be written as:
$$-D = \prod_{i=1}^t q_i^*$$
where all $q_i^*$'s are distinct and
$q_i^*$ is either $-4$ or $\pm 8$, or $q_i^* = (-1/q_i) q_i$
for any prime $q_i$. The number of genera is $g(-D) = 2^{t-1}$ and
Gauss proved that this number divides the class number $h(-D)$ of
$\K=\Q(\sqrt{-D})$. Moreover, Siegel proved that $h =
O(D^{1/2+\varepsilon})$ asymtotically.

The rational prime $p$ is the norm of an integer in $\K$, or
equivalently, $4 p = U^2+D V^2$ in rational integers $U$ and $V$ if
and only if the ideal $(p)$ splits completely in the Hilbert class
field of $\K$, denoted $\KH$, an extension of degree $h(-D)$ of $\K$.
The probability that a prime $p$ splits in $\K$
is $1/(2 h(-D))$. 

Using Gauss's theory of genera of forms, it is known that
if $\legendre{q_i^*}{p} = 1$ for all $i$ (equivalently, $(p)$ splits
in the genus field of $\K$), then
the probability of $(p)$ splitting in $\KH$ is $g(-D)/h(-D)$.

%%%%% S
\section{The basic ECPP algorithm}

We present a rough sketch of the ECPP algorithm, enough for us to
estimate its complexity. We do not insist on what happens if one of
the steps fails, revealing the compositeness of $N$. More details can
be found in \cite{AtMo93b}.

%%%%%%%%%% SS
\subsection{Elliptic curves over $\ZN{N}$}

For us, an elliptic curve $E$ modulo $N$ will have an equation $Y^2
\equiv X^3+a X+b$ with $\gcd(4a^3+27 b^2, N)=1$ and we will use the
set of points $E(\ZN{N})$ defined as:
$$E(\ZN{N}) = \{(x:y:z) \in \mathbb{P}^2(\ZN{N}), y^2 z \equiv x^3 + a
x z^2 + b z^3\} \cup \{O_E = (0:1:0)\}$$
which ressembles the definition of an actual elliptic curve if $N$ is
prime, $\mathbb{P}^2(\ZN{N})$ being the projective plane over $\ZN{N}$. The
important point here is that if $p$ is a divisor of $N$, we can reduce
the curve $E$ and a point $P$ on it via a reduction modulo $p$ of each
integer, yielding a point $P_p$ on $E_p$. Moreover, we can define an
operation on $E(\ZN{N})$, called {\em pseudo-addition}, that adds two
``points'' $P$ and $Q$ with the usual chord-and-tangent law. This
operation either yields a point $R$ or a divisor of $N$ if any is
encountered when dividing. If $R$ exists, then it has the property
that $R_p$ is the sum of $P_p$ and $Q_p$ on $E_p$ for all prime
factors $p$ of $N$. Note also that $O_E$
reduces to the ordinary point at infinity on $E_p$.

We will need to exponentiate points in $E$. This is best defined using
the division polynomials (see for instance \cite{Ayad92} for a lot of
properties on these). Remember that over a field $K$ there exist
polynomials $\phi_m(X, Y)$, $\psi_m(X, Y)$, $\omega_m(X, Y)$ such that
\begin{equation}\label{div-exp}
[m] P = \underbrace{P+\cdots+P}_{m\text{ times}} = [m] (X, Y) =
\left(\phi_m(X, Y)\psi_m(X, Y): \omega_m(X, Y): \psi_m^3(X,
Y)\right).
\end{equation}
All these polynomials can be computed via recurrence formulas and
there is a $O(\log m)$ algorithm for this task (a variant of the usual
binary method for exponentiating). 

We will take (\ref{div-exp}) for the definition of $[m] P$ over
$\ZN{N}$. We note here that if $\psi_m(X, Y) = 0$, then $[m] P$ is
equivalent to the point $O_E$.

For the sake of presenting the algorithm in a simplified setting, we
prove (compare \cite{GoKi99}):
\begin{proposition}\label{prop-prim}
Let $N'$ a prime satisfying $(\sqrt{N}-1)^2 \leq 2 N' \leq
(\sqrt{N}+1)^2$. Suppose that $E(\ZN{N})$ is a curve over $\ZN{N}$,
that $P=(x:y:1)$ is such that $\gcd(y, N)=1$, 
$\psi_{2N'}(x, y) = 0$ but $\gcd(\psi_{N'}(x, y), N)=1$. Then $N$ is prime.
\end{proposition}

\noindent
{\em Proof:} suppose that $N$ is composite and that $p \leq \sqrt{N}$
is one of its prime factors. Let us look at what happens modulo $p$. By
construction, $P_p$ is not a 2-torsion point on $E_p$. Since
$\psi_{N'}(x, y)$ is invertible modulo $p$, then $[N'] (P_p) \neq
O_{E_p}$ and therefore $P_p$ is of order $2N'$ modulo $p$. This is
impossible, since $2 N' \geq (\sqrt{N}-1)^2 \geq (p-1)^2 >
(\sqrt{p}-1)^2 \geq \#E_p$ by Hasse's theorem. $\Box$

%%%%%%%%%% SS
\subsection{Presentation of the algorithm}

We want to prove that $N$ is prime. The algorithm runs as follows:

\medskip
\noindent [Step 1.] Repeat the following: Find an imaginary quadratic
field $\K = \Q(\sqrt{-D})$ of discriminant $-D$, $D > 0$, such that
\begin{equation}\label{eqnorm}
4 N = U^2 + D V^2
\end{equation}
in rational integers $U$ and $V$. 
For all solutions $U$ of (\ref{eqnorm}), compute
$m = N+1-U$; if one of these numbers is twice a probable prime $N'$,
go to Step 2.

\medskip
\noindent [Step 2.] Build an elliptic curve $\overline{E}$ over
$\overline{\Q}$ having complex multiplication by the ring of integers $\OK$
of $\K$.

\medskip
\noindent [Step 3.] Reduce $\overline{E}$ modulo $N$ to get a curve
$E$. 

\noindent [Step 4.] Find $P = (x:y:1)$, $\gcd(y, N)=1$ on $E$
such that $\psi_{2N'}(x, y)=0$, but $\gcd(\psi_{N'}(x, y), N)=1$.
If this cannot be done, then $N$ is composite, otherwise, it is prime
by Proposition \ref{prop-prim}.

\medskip
\noindent [Step 5.] Set $N = N'$ and go back to Step 1.

%%%%%%%%%% SS
\subsection{Analyzing ECPP}

We will now analyze all steps of the above algorithm and give complexity
estimates using the parameter $L = \log N$. One basic unit
of time will be the time needed to multiply two integers of size $L$,
namely $O(L^{1+\mu})$, where $0\leq\mu\leq 1$ ($\mu=1$ for ordinary
multiplication, or $\epsilon > 0$ for any fast multiplication
method).

Clearly, we need $\log N$ steps for proving the primality of $N$. 
We consider all steps, one at a time, easier steps first.

%%%%%%%%%% SS
\subsection{Analysis of Step 4.}

Finding a point $P$ can be done by a simple algorithm that looks for
the smallest $x$ such that $x^3+a x+b$ is a square modulo $p$ and then
extracting a squareroot modulo $p$, for a cost of $O((\log
N)^{2+\mu})$. Note that we can do without this with the trick
described in \cite[\S 8.6.3]{AtMo93b}, though we do not need this at
this point. 

Computing $\psi_{N'}(x, y)$ costs $O((\log N)^{2+\mu})$, and we need
$O(1)$ points on average, so this steps amounts for $O((\log N)^{2+\mu})$.

%%%%%%%%%% SS
\subsection{Analyzing Step 2}

The original version is to realize $\KH/\K$ via the
computation of the minimal polynomial $H_D(X)$ of the special values
of the classical $j$-invariant at quadratic integers. More precisely,
we can view the class group $Cl(-D)$ of $\K$ as a set of primitive reduced
quadratic forms of discriminant $-D$. If $(A, B, C)$ is such a form,
with $B^2 - 4 A C = -D$, then
$$H_D(X) = \prod_{(A, B, C) \in Cl(-D)} \left(X -
j((-B+\sqrt{-D})/(2A))\right).$$
In \cite{EnMo02}, it is argued that the height of this
polynomial is well approximated by the quantity:
\[
\pi\sqrt{D} \sum_{[A, B, C] \in Cl(-D)} \frac{1}{A},
\]
which can be shown to be $O((\log h)^2)$.

Evaluating the roots of $H_D(X)$ and building this polynomial can be
done in $\tilde{O}(h^2)$ operations (see \cite{Enge04}). Note that
this step does not require computations modulo $N$.

Alternatively, we could use the method of \cite{CoHe02,BrSt04} for
computing the class polynomial and get a proven running time of
$\tilde{O}(h^2)$, but assuming GRH.

%%%%%%%%%% SS
\subsection{Analyzing Step 3}

Reducing $E$ modulo $N$ is done by finding a root of
$H_D(X)$ modulo $N$. This can be done with
the Cantor-Zassenhaus algorithm (see \cite{GaGe99} for
instance). Briefly, we split recursively $H_D(X)$ by
computing $\gcd((X+a)^{(N-1)/2}-1, H_D(X))\bmod N$ for random $a$'s.

Computing $(X+a)^{(N-1)/2} \bmod (N, H_D(X))$ costs
$O((\log N) {\sf M}(N, h)) = O(L {\sf M}(N, h))$ where ${\sf M}(N, d)$
is the time needed to multiply two degree $d$ polynomials modulo
$N$. A gcd of two degree $d$ polynomials costs ${\sf M}(N, d)\log d$ (see
\cite[Ch. 11]{GaGe99}). % pp. 304, Coro 11.6
The total splitting requires $\log h$ steps, but the overall cost
is dominated by the first one, hence yields a time:
$$O({\sf M}(N, h) \max(L, \log h)).$$
We can assume that ${\sf M}(N, d) = O(d^{1+\nu} L^{1+\mu})$ where
again $0\leq \nu\leq 1$.

%%%%%%%%%% SS
\subsection{Analysis of Step 1.}

This is the crucial step that will give us the clue to the complexity.
Given $D$, testing whether (\ref{eqnorm}) is satisfied involves the
reduction of the ideal $(N, \frac{r-\sqrt{-D}}{2})$ that lies above
$(N)$ in $\K$, where $r^2 \equiv -D \bmod (4N)$ (if $N$ is
prime...). This requires the computation of $\sqrt{-D}\bmod N$, using
for instance the Tonelli-Shanks algorithm, for the cost of one modular
exponentiation, i.e., a $O(L^{2+\mu})$ time. Then it proceeds with a
half gcd like computation, for a cost of $O(L^{1+\mu})$ (see also
section \ref{ssct:Cornacchia} below).

In the event that equation (\ref{eqnorm}) is solvable, then we need
check that $m = 2 N'$ and test $N'$ for primality, which costs again some
$O(L^{2+\mu})$.

The heuristic probability of $m$ being of the given form is
$O(1/L)$. Though quite realistic, it is impossible to prove, given the
current state of the art in analytical number theory.
Using this heuristics, we expect to need $O(L)$ splitting $D$'s.
Let us take all discriminants less than $D_{\max}$. They have class
number close to $h(-D_{\max}) = O(\sqrt{D_{\max}})$ and there are
$O(D_{\max})$ of them. We see that if $L =
O(D_{\max})/\sqrt{D_{\max}}$, then among these discriminants,
one will lead to a useful $m$. We conclude that $D_{\max} = O(L^2)$
should suffice.

Turning to complexity, the cost of Step 1 is then that of $O(L^2)$
solving of (\ref{eqnorm}), followed by $O(L)$ probable primality tests:
$$O(L^{2} (\underbrace{L^{2+\mu}}_{\sqrt{-D} \bmod N} +
\underbrace{L^{1+\mu}}_{\text{reduction}}))
+ O(L \cdot \underbrace{L^{2+\mu}}_{\text{probable primality}}).$$
which is dominated by the first cost, namely $O(L^{4+\mu})$.

%%%%%%%%%% SS
\subsection{Adding everything together}

Taking $D = O(L^2)$ readily implies $h = O(L)$, so that the cost of
Step 2 is $\tilde{O}(L^2)$, and that of Step 3 is $O((\log L)
L^{3+\mu+\nu})$, which dominates Step 4. All in all, we get that ECPP
has heuristic complexity $O(L^{4+\mu})$ for one step, and therefore
$O(L^{5+\mu})$ in totality.

%%%%%%%%%% SS
\subsection{Remark}

In practice, the dominant term of the complexity of Step 1 is $O(n_D
L^{2+\mu})$ where $n_D$ is the number of $D$'s for which we try solve
equation (\ref{eqnorm}). Depending on implementation parameters and
real size of $N$, this number $n_D$ can be quite small. This gives a
very small {\em apparent complexity} to ECPP, somewhere in between
$L^{3}$ and $L^4$ and explains why ECPP seems so fast in practice (see
for instance \cite{Galway04}).

%%%%% S
\section{The fast version of ECPP}

%%%%%%%%%% SS
\subsection{Presentation}

When dealing with large numbers, all the time is spent in the finding
of $D$, which means that a lot of squareroots modulo $N$ must be
computed. A first way to reduce the computations, alluded to in
\cite[\S 8.4.3]{AtMo93b}, is to accumulate squareroots, and reuse
them, at the cost of some multiplications. For instance, if one has
$\sqrt{-3}$ and $\sqrt{5} = \sqrt{-20}/\sqrt{-4}$, then we can build
$\sqrt{-15}$, etc.

A better way that leads to the fast version consists in computing
a basis of small squareroots and build discriminants from this
basis. Looking at the analysis carried out above, we see that we
need $O(L^2)$ discriminants to find a good one. The basic version
finds them by using all discriminants that are of size $O(L^2)$. As
opposed to this, one can build those discriminants as $-D = (-p) (q)$,
where $p$ and $q$ are taken from a pool of size $O(L)$ primes.

More formally, we replace Step 1. by Step 1'. as follows:

\medskip
\noindent [Step 1'.] 

1.1. Find the $r = O(L)$ smallest primes $q^*$ such that
   $\legendre{q^*}{N} = 1$, yielding $\mathcal{Q} = \{q_1^*, q_2^*,
   \ldots, q_r^* \}$.

1.2. Compute all $\sqrt{q^*}\bmod N$ for $q^*\in\mathcal{Q}$.

1.3. For all pairs $(q_{i_1}^*, q_{i_2}^*)$ of $\mathcal{Q}$ for which
   $q_{i_1}^* q_{i_2}^* = -D < 0$, try to solve equation
   (\ref{eqnorm}).

\medskip
The cost of this new Step 1 is that of computing $r = O(L)$
squareroots modulo $N$, for a cost of $O(L \cdot L^{2+\mu})$. Then, we
still have $O(L^2)$ reductions. The new overall cost of this phase
decreases now to:
$$O(L\cdot \underbrace{L^{2+\mu}}_{\text{squareroots}})
+ O(L^2 \cdot \underbrace{L^{1+\mu}}_{\text{reduction}})
+ O(L\cdot \underbrace{L^{2+\mu}}_{\text{probable primality}})$$
which yields namely $O(L^{3+\mu})$.
Note here how the complexity decomposes as $3 =1+2$ or $2+1$ depending
on the sub-algorithms.

Putting everything together, we end up with a total cost of
$O(L^{4+\mu})$ for this variant of ECPP.

%%%%%%%%%%
\subsection{Remarks}

%%%%%%%%%%%%%%% SSS
\subsubsection{Complexity issues}

We can slightly optimize the preceding argument, by using all subsets
of $\mathcal{Q}$ and not only pairs of elements. This would call for
$r = O(\log\log N)$, since then $2^r = L^2$ could be reached. Though
useful in practice, this phase no longer dominates the cost of the
algorithm. 

Moreover, we can see that several phases of {\sc fastECPP} have cost
$\tilde{O}(L^3)$, which means that we would have to fight hard to
decrease the overall complexity below $\tilde{O}(L^{4})$.

%%%%%%%%%%%%%%% SSS
\subsubsection{A note on discriminants}

Note that we use fundamental discriminants only, as non fundamental
discriminants lead to curves that do not bring anything new compared
to fundamental ones. Indeed, if $\mathcal{D} = f^2 D$, with $D$
fundamental, then there is a curve having CM by the order of
discriminant $\mathcal{D}$. Writing $4 N = U^2 + D f^2 V^2$, its
cardinality is $N+1-U$, the same as the corresponding curve associated
to $D$.

%%%%%%%%%%%%%%% SSS
\subsubsection{A note on class numbers}

As soon as we use composite discriminants $-D$ of the form $q_{i_1}^*
q_{i_2}^*$, Gauss's theorem tells us that the class number $h(-D)$ is
even. This could bias our estimation, but we conjecture that the
effect is not important.

%%%%% S
\section{Implementation}

%%%%%%%%%% SS
\subsection{Computing class numbers}

In order to make the search for $D\in\mathcal{D}$ efficient, it is
better to control the class number beforehand. Tables can be made, but
for larger computations, we need a fast way to compute
$h(-D)$. Subexponential methods exist, assuming the Generalized
Riemann Hypothesis. From a practical point of view, our $D$'s are
of medium size. Enumerating all forms costs $O(h^2)$ with a small
constant, and Shanks's baby-steps/giant-steps algorithm costs
$O(\sqrt{h})$ but with a large constant. It is better here to use
the explicit formula of Louboutin \cite{Louboutin02} that yields a
practical method in $O(h)$ with a very small constant.

%%%%%%%%%% SS
\subsection{An improved Cornacchia algorithm}\label{ssct:Cornacchia}

Step 1 needs squareroots to be computed, and some half gcd to be
performed. Briefly, Cornacchia's algorithm runs as follows (see
\cite{Nitaj95}):

\medskip
\noindent {\bf procedure} {\sc Cornacchia}($d$, $p$, $t$)

\{$t$ is such that $t^2 \equiv -d \bmod p$, $p/2 < t < p$ \}

\ \ a) $r_{-2} = p$, $r_{-1} = t$~; $w_{-2} = 0$, $w_{-1} = 1$;

\ \ b) {\bf for} $i\geq 0$ {\bf while} $r_{i-1} > \sqrt{p}$ {\bf do}

\ \ \ \ \ \ \ $r_{i-2} = a_i r_{i-1} + r_i$, $0 \leq r_i < r_{i-1}$~;

\ \ \ \ \ \ \ $w_i = w_{i-2}+a_i w_{i-1}\quad (*)$~;

\ \ c) {\bf if} $r_{i-1}^2 + d w_{i-1}^2 = p$ {\bf then}
return $(r_{i-1}, w_{i-1})$ {\bf else} return $\emptyset$.

\medskip
We end the {\bf for} loop once we get $r_{i-1} \leq \sqrt{p} <
r_{i-2}$. As is well known, the $a_i$'s are quite small and we
can guess their size by monitoring the number of bits of the $r_i$'s,
thus limiting the number of long divisions. One can use a fast variant
for this half gcd if needed, in a way reminiscent of Knuth. 

Moreover, from the theory, we know that this algorithms almost always
returns that the empty set in step 2c), since the probability of
success if $1/(2h(-d))$. Therefore, when $h$ is large, we can dispense of the
multiprecision computations in equation $(*)$. We replace it by single
precision computations:
$$w_i = w_{i-2}+a_i w_{i-1} \bmod 2^{32}$$
and at the end, we test whether $r_{i-1}^2 + d w_{i-1}^2 = p\bmod
2^{32}$. If this is the case, then we redo the computation of the
$w_i$'s and check again.

%%%%%%%%%% SS
\subsection{Factoring $m$}

Critical parameters are that related to the factorization of $m$,
since in practice we try to factor $m$ to get it of the form $c N'$
for some $B$-smooth number $c$.

As shown in \cite{FrKlMoWi04}, the number of probable prime tests we
will have to perform is $t = O((\log N)/(\log B))$ and we will end up with
$N'$ such that $N/N' \approx \log B$.

For small numbers, we can factor lots of $m$ doing the following.
In a first step, we compute
$$r_i = (N+1) \bmod p_i$$
for all $p_i\leq B$'s, which costs $\pi(B) (\log N)^{1+\varepsilon}$,
where $\pi(B) = O(B/\log B)$ is the number of primes below $B$ and the
other term being the time needed to divide a multi-digit number by a
single digit number.

Then, sieving both $m = N+1-U$ and $m' = N+1+U$ is
done by computing $u_i = U \bmod p_i$ and comparing it to $\pm
r_i$ for primes $p_i$ such that $(-D/p_i)\neq -1$. See
\cite{AtMo93b,Morain98a} for more details and tricks. 

The cost of this algorithm for $t$ values of $m$ is
$$O(\pi(B) (\log N)^{1+\varepsilon}) + O(t \cdot \pi(B) (\log
N)^{1+\varepsilon})$$
where the second term is that for computing $U\bmod p_i$, which is
slightly half that of $(N+1)\bmod p_i$, since $U$ is
$O(\sqrt{N})$. Since we need to perform also $t$ probable prime tests
(say, a plain Fermat one), then the cost is
$$O(t B L) + O(t L^{2+\mu}) = O(B L^2) + O(L^{3+\mu})$$
and therefore the optimal value for $B$ is $B = O(L)$. 

For larger numbers, it is better to use the stripping factor algorithm in
\cite{FrKlMoWi04}, for a cost of $O(B (\log B)^2)$, the optimal value
of $B$ being $B = O((\log N)^3)$.

%%%%%%%%%%%%%%% SSS
\subsubsection{Remark}

Suppose now that we have found $N'$ and that $m = N+1-U = c N'$. Then
we will have to compute 
$$r_i' = (N'+1) \bmod p_i$$
which may be computed as:
$$r_i' = (r_i-u_i)/c+1 \bmod p_i.$$
Computing the right hand side is faster, since $c$ is ordinarily small
compared to $N'$.

%%%%%%%%%%%%%%% SSS
\subsubsection{Using an early abort strategy}\label{delta}

This idea is presented in \cite{FrKlMoWi04}.
We would like to go down as fast as possible. So why not impose $N/N'$
greater than some given bound? Candidates $N'$ need be tested for
probable primality only if this bound is met. From what has been
written above, we can insist on $N/N' \approx \log B$. In practice, we
used a bound $\delta$ and used $N/N' \geq 2^{\delta}$.

%%%%%%%%%%%%%%% SSS
\subsubsection{Using new invariants}

Proving larger and larger numbers forces us to use larger and larger
$D$'s, leading to larger and larger polynomials $H_D$.
For this to be doable, new invariants had to be used, so as to
minimize the size of the minimal polynomials. This task was done using
Schertz's formulation of Shimura's reciprocity law \cite{Schertz02},
with the invariants of \cite{EnSc03} as demonstrated in
\cite{EnMo02} (alternatively see
\cite{GeSt98,Gee99}). Note that replacing $j$ by other functions does
not the change the complexity of the algorithm, though it is crucial
in practice.

%%%%%%%%%%%%%%% SSS
\subsubsection{Step 3 in practice}

We already noted that this step is the theoretically dominating one in
{\sc fastECPP}, with a cost of $O((\log L) L^{3+\mu+\nu})$. In practice,
even for small values of $h$, we can assume $\nu \approx 0$ (using for
instance the algorithm of \cite{Mihailescu04} for polynomial multiplication).

Galois theory comes in handy for reducing the $\log L$ term to a
$\log\log L$ one, if we insist on $h$ being smooth. Then, we replace
the time needed to factor a degree $h$ polynomial by a list of smaller
ones, the largest prime factor of $h$ being $\log h$. We already used
that in ECPP, using \cite{HaMo01,EnMo03}. Typical values of $h$ are
now routinely in the $10000$ zone.

It could be argued that keeping only smooth class numbers is too
restrictive. Note however, that class numbers tend to be smoother than
ordinary numbers \cite{CoLe84b}.

\ifdraft
\else
Atkin's trick?
\fi

%%%%%%%%%%%%%%% SSS
\subsubsection{Improving the program}

The new implementation uses {\sc
GMP}\footnote{http://www.swox.com/gmp/} for the basic arithmetic,
which enables one to use \verb+mpfr+ \cite{HaLeZial02} and \verb+mpc+
\cite{EnZi02}, thus leading to a
complete program that can compute polynomial $H_D$'s on the fly,
contrary to the author's implementation of ECPP, prior to version
11.0.5. This turned out to be the key for the new-born program to compete
with the old one.

%%%%%%%%%% SS
\subsection{{\sc fastECPP}}

We give here the expanded algorithm corresponding to step 1'. Using a
smoothness bound $B$, we need approximately $t = \exp(-\gamma) \log
N/\log B$ values of $m$ and therefore roughly $t/2$ discriminants. The
probability that $D$ is a splitting discriminant is
$g(-D)/h(-D)$. Therefore we build discriminants until
$$\sum_D g(-D)/h(-D) \approx t/2.$$
One way of building these discriminants is the following: we let $r$
increase and build all or some of the subsets of $\{q_1^*, \ldots,
q_r^*\}$ until the expected number of $D$'s is reached. After this, we
sort the discriminants with respect to $(h(-D)/g(-D), h(-D), D)$ and
treat them in this order.

%%%%% S
\section{Benchmarks}

First of all, it should be noted that ECPP is not a well defined
algorithm, as long as one does not give the list of discriminants that
are used, or the principles that generate them.

Since the first phase of ECPP requires a tree search,
testing on a single number does not reveal too much. Averaging on more
than 20 numbers is a good idea.

Our current implementation uses {\sc
GMP}\footnote{http://www.swox.com/gmp/} for the basic arithmetic,
which enables one to use \verb+mpfr+ \cite{HaLeZial02} and \verb+mpc+
\cite{EnZi02}, thus leading to a
complete program that can compute polynomial $H_D$'s on the fly,
contrary to the author's implementation of ECPP, prior to version
11.0.5. This turned out to be the key for the new-born program to compete
with the old one.

We give below some timings obtained with this implementation, after a
lot of trials. We used
as prime candidates the first twenty primes of 1000, 1500, and 2000
decimal digits. Critical parameters are as follows: we used
$D \leq 10^7$, $h \leq 1000$, $\delta = 12$ (see section
\ref{delta}). For 1000 and
1500 decimal digits, we limited the largest prime factor of $h$ to be
$\leq 30$ and for 2000 dd, it was put to $100$. This parameter has an
influence in Step 3. For the extraction of small prime factors (used
in the algorithm described in \cite{FrKlMoWi04} and denoted EXTRACT in
the sequel), we used $B = 8\cdot 10^6$, $10^7$, $3\cdot 10^7$ for the
three respective sizes.

SQRT refers to the computation of the $\sqrt{q_i^*}$, CORN to
Cornacchia, PRP to probable primality tests; HD is the time
for computing polynomials $H_D$ using the techniques described in
\cite{EnMo02}, jmod the time to solve it modulo $p$; then 1st refers
to the building phase (step 1), 2nd to the other ones; total is the
total time, check the time to verify the certificate. Follow some
data concering $D$, $h$ and the size of the certificates (in kbytes).
All timings are cumulated CPU time on an AMD Athlon 64 3400+ running
at 2.4GHz.

\begin{center}
\begin{table}[hbt]
\begin{tabular}{|l|r|r|r|r|}\hline
 & min & max & avg & std \\\hline
SQRT   &        19 &    34 &    25 &    3\\
CORN   &        10 &    24 &    17 &    4\\
EXTRACT   & 60 &    84 &    74 &    5\\
PRP   & 74 &    124 &   102 &   14\\
\hline
HD   &  0 &     7 &     2 &     2\\
jmod   &        42 &    99 &    61 &    11\\
\hline
1st   & 178 &   276 &   234 &   27\\
2nd   & 79 &    136 &   99 &    12\\
total   &       260 &   387 &   334 &   34\\
\hline
check   &       18 &    22 &    20 &    0\\
\hline
nsteps   &      124 &   156 &   143 &   7\\
certif   &      396 &   456 &   435 &   13\\
D & & 8740947 & 120639 & 608050 \\
h & & 1000 & 31 & 87 \\
\hline
\end{tabular}
\caption{1000 decimal digits}
\end{table}
\end{center}

\begin{center}
\begin{table}[hbt]
\begin{tabular}{|l|r|r|r|r|}\hline
 & min & max & avg & std \\\hline
SQRT   &        114 &   427 &   171 &   65\\
CORN   &        59 &    140 &   95 &    21\\
EXTRACT   & 195 &   282 &   230 &   20\\
PRP   & 472 &   903 &   664 &   99\\
\hline
HD   &  5 &     13 &    9 &     2\\
jmod   &        219 &   471 &   334 &   60\\
\hline
1st   & 868 &   1590 &  1192 &  185\\
2nd   & 368 &   649 &   508 &   70\\
total   &       1322 &  2240 &  1701 &  230\\
\hline
check   &       71 &    94 &    85 &    5\\
\hline
nsteps   &      183 &   209 &   198 &   7\\
certif   &      796 &   968 &   897 &   40\\
D & & 9644776 & 201015 & 848112 \\
h & & 972 & 46 & 111 \\
\hline
\end{tabular}
\caption{1500 decimal digits}
\end{table}
\end{center}

\begin{center}
\begin{table}[hbt]
\begin{tabular}{|l|r|r|r|r|}\hline
 & min & max & avg & std \\\hline
SQRT   &        384 &   820 &   516 &   120\\
CORN   &        181 &   390 &   260 &   55\\
EXTRACT   & 600 &   853 &   713 &   67\\
PRP   & 1761 &  2879 &  2227 &  306\\
\hline
HD   &  6 &     27 &    16 &    5\\
jmod   &        969 &   1539 &  1255 &  188\\
\hline
1st   & 2974 &  4888 &  3778 &  528\\
2nd   & 1398 &  2120 &  1777 &  221\\
total   &       4494 &  6795 &  5557 &  711\\
\hline
check   &       213 &   261 &   238 &   13\\
\hline
nsteps   &      236 &   262 &   248 &   7\\
certif   &      1420 &  1644 &   1539 & 64 \\
D & & 9760387 & 285217 & 1026529 \\
h & & 1000 & 63 & 130 \\
\hline
\end{tabular}
\caption{2000 decimal digits}
\end{table}
\end{center}

Looking at the average total time, we see that it follows very closely
the $O((\log N)^4)$ prediction. Note also that the dominant time is
that of the PRP tests, and that all phases have time close to what was
predicted.

%%%%% S
\section{Conclusions}

We have described in greater details the fast version of ECPP. We have
demonstrated its efficiency. As for ECPP, it is obvious that the
computations can be distributed over a network of computers. We refer
the reader to \cite{FrKlMoWi04} for more details. Note that the
current record of 15041 decimal digits (with the number
$4405^{2638}+2638^{4405}$ see transaction in the NMBRTHRY mailing list),
was settled using this approach. Many more numbers were proven prime
using either the monoprocessor version or the distributed one, most of
them from the tables of numbers of the form
$x^y+y^x$ made by P.~Leyland\footnote{{\tt
http://www.leyland.vispa.com/numth/primes/xyyx.htm}}.

Cheng \cite{Cheng03} has suggested to use ECPP to help his improvement
of the AKS algorithm, forcing $m = c N'$ to have $N'-1$ divisible by a
given prime large prime of size $O((\log N)^2)$. The same idea can be
used to speed up the Jacobi sums algorithm, and this will be detailed
elsewhere.

\medskip
\noindent {\bf Acknowledgments.} The author wants to thank N.~Bourbaki
for making him dive once again in the field of primality proving and
D.~Bernstein for stimulating emails on the existence and analysis of
{\sc fastECPP}. My co-authors of \cite{FrKlMoWi04} were a source of
stimulation through their records.
Thanks also to P.~Gaudry for never ending discussions on
how close we are to infinity, as far as fast algorithms are
concerned. D.~Stehl\'e and P.~Zimmermann for useful discussions around
Cornacchia and fast sieving. Thanks to A.~Enge for his help in
improving the exposition, and to D.~Bernardi for his remarks that
helped clarify the exposition.

\iffalse
\bibliographystyle{plain}
\bibliography{morain}

\begin{thebibliography}{10}

\bibitem{AdPoRu83}
L.~M. Adleman, C.~Pomerance, and R.~S. Rumely.
\newblock On distinguishing prime numbers from composite numbers.
\newblock {\em Ann. of Math. (2)}, 117:173--206, 1983.

\bibitem{AgKaSa02}
M.~Agrawal, N.~Kayal, and N.~Saxena.
\newblock {PRIMES} is in {P}.
\newblock \Preprint; available at {\tt
  http://www.cse.iitk.ac.in/primality.pdf}, August 2002.

\bibitem{AtMo93b}
A.~O.~L. Atkin and F.~Morain.
\newblock Elliptic curves and primality proving.
\newblock {\em Math. Comp.}, 61(203):29--68, July 1993.

\bibitem{Ayad92}
M.~Ayad.
\newblock Points ${S}$-entiers des courbes elliptiques.
\newblock {\em Manuscripta Math.}, 76(3-4):305--324, 1992.

\bibitem{Bernstein03a}
D.~Bernstein.
\newblock Proving primality after {A}grawal-{K}ayal-{S}axena.
\newblock {\tt http://cr.yp.to/papers/aks.ps}, January 2003.

\bibitem{Bernstein03}
D.~Bernstein.
\newblock Proving primality in essentially quartic expected time.
\newblock {\tt http://cr.yp.to/papers/quartic.ps}, January 2003.

\bibitem{Berrizbeitia02}
P.~Berrizbeitia.
\newblock Sharpening "{P}rimes is in {P}" for a large family of numbers.
\newblock {\tt http://arxiv.org/abs/math.NT/0211334}, November 2002.

\bibitem{BrSt04}
R.~Br{\"o}ker and P.~Stevenhagen.
\newblock Elliptic curves with a given number of points.
\newblock In D.~Buell, editor, {\em Algorithmic Number Theory}, volume 3076 of
  {\em Lecture Notes in Comput. Sci.}, pages 117--131. Springer-Verlag, 2004.
\newblock 6th International Symposium, ANTS-VI, Burlington, VT, USA, June 2004,
  Proceedings.

\bibitem{Cheng03}
Q.~Cheng.
\newblock Primality proving via one round in {ECPP} and one iteration in {AKS}.
\newblock In D.~Boneh, editor, {\em Advances in Cryptology -- CRYPTO 2003},
  volume 2729 of {\em Lecture Notes in Comput. Sci.}, pages 338--348. Springer
  Verlag, 2003.

\bibitem{CoLe84b}
H.~Cohen and H.~W. {Lenstra, Jr.}
\newblock Heuristics on class groups of number fields.
\newblock In H.~Jager, editor, {\em Number Theory, Noordwijkerhout 1983},
  volume 1068 of {\em Lecture Notes in Math.}, pages 33--62. Springer-Verlag,
  1984.
\newblock Proc. of the Journ{\'e}es Arithm{\'e}tiques 1983, July 11--15.

\bibitem{CoLe84}
H.~Cohen and H.~W. {Lenstra, Jr.}
\newblock Primality testing and {J}acobi sums.
\newblock {\em Math. Comp.}, 42(165):297--330, 1984.

\bibitem{CoHe02}
J.-M. Couveignes and T.~Henocq.
\newblock Action of modular correspondences around {CM} points.
\newblock In C.~Fieker and D.~R. Kohel, editors, {\em Algorihmic Number
  Theory}, volume 2369 of {\em Lecture Notes in Comput. Sci.}, pages 234--243.
  Springer-Verlag, 2002.
\newblock 5th International Symposium, ANTS-V, Sydney, Australia, July 2002,
  Proceedings.

\bibitem{Cox89}
D.~A. Cox.
\newblock {\em Primes of the form $x^2+n y^2$}.
\newblock John Wiley \& Sons, 1989.

\bibitem{CrPo00}
R.~Crandall and C.~Pomerance.
\newblock {\em Prime numbers -- A Computational Perspective}.
\newblock Springer Verlag, 2000.

\bibitem{Enge04}
A.~Enge.
\newblock The complexity of class polynomial computations via floating point
  approximations.
\newblock Preprint, February 2004.

\bibitem{EnMo02}
A.~Enge and F.~Morain.
\newblock Comparing invariants for class fields of imaginary quadratic fields.
\newblock In C.~Fieker and D.~R. Kohel, editors, {\em Algorithmic Number
  Theory}, volume 2369 of {\em Lecture Notes in Comput. Sci.}, pages 252--266.
  Springer-Verlag, 2002.
\newblock 5th International Symposium, ANTS-V, Sydney, Australia, July 2002,
  Proceedings.

\bibitem{EnMo03}
A.~Enge and F.~Morain.
\newblock Fast decomposition of polynomials with known {G}alois group.
\newblock In M.~Fossorier, T.~H{\o}holdt, and A.~Poli, editors, {\em Applied
  Algebra, Algebraic Algorithms and Error-Correcting Codes}, volume 2643 of
  {\em Lecture Notes in Comput. Sci.}, pages 254--264. Springer-Verlag, 2003.
\newblock 15th International Symposium, AAECC-15, Toulouse, France, May 2003,
  Proceedings.

\bibitem{EnSc03}
A.~Enge and R.~Schertz.
\newblock Modular curves of composite level.
\newblock Soumis, 2003.

\bibitem{EnZi02}
A.~Enge and P.~Zimmermann.
\newblock {\tt mpc} --- a library for multiprecision complex arithmetic with
  exact rounding, 2002.
\newblock Version 0.4.1, available from
  http://www.lix.polytechnique.fr/Labo/Andreas.Enge.

\bibitem{FrKlMoWi04}
J.~Franke, T.~Kleinjung, F.~Morain, and T.~Wirth.
\newblock Proving the primality of very large numbers with fastecpp.
\newblock In D.~Buell, editor, {\em Algorithmic Number Theory}, volume 3076 of
  {\em Lecture Notes in Comput. Sci.}, pages 194--207. Springer-Verlag, 2004.
\newblock 6th International Symposium, ANTS-VI, Burlington, VT, USA, June 2004,
  Proceedings.

\bibitem{Galway04}
W.~F. Galway.
\newblock {\em Analytic computation of the prime-counting function}.
\newblock PhD thesis, University of Urbana-Champaign, 2004.
\newblock http://www.math.uiuc.edu/\~{}galway/PhD\_Thesis/.

\bibitem{GaGe99}
{J. von zur} Gathen and J.~Gerhard.
\newblock {\em Modern Computer Algebra}.
\newblock Cambridge University Press, 1999.

\bibitem{Gee99}
A.~Gee.
\newblock Class invariants by {S}himura's reciprocity law.
\newblock {\em J. Th\'eor. Nombres Bordeaux}, 11:45--72, 1999.

\bibitem{GeSt98}
A.~Gee and P.~Stevenhagen.
\newblock Generating class fields using {S}himura reciprocity.
\newblock In J.~P. Buhler, editor, {\em Algorithmic Number Theory}, volume 1423
  of {\em Lecture Notes in Comput. Sci.}, pages 441--453. Springer-Verlag,
  1998.
\newblock Third International Symposium, ANTS-III, Portland, Oregon, june 1998,
  Proceedings.

\bibitem{GoKi99}
S.~Goldwasser and J.~Kilian.
\newblock Primality testing using elliptic curves.
\newblock {\em Journal of the ACM}, 46(4):450--472, July 1999.

\bibitem{HaLeZial02}
G.~Hanrot, V.~Lef\`evre, and P.~{Zimmermann et. al.}
\newblock {\tt mpfr} --- a library for multiple-precision floating-point
  computations with exact rounding, 2002.
\newblock Version contained in {\sc gmp}. Available from http://www.mpfr.org.

\bibitem{HaMo01}
G.~Hanrot and F.~Morain.
\newblock Solvability by radicals from an algorithmic point of view.
\newblock In B.~Mourrain, editor, {\em Symbolic and algebraic computation},
  pages 175--182. ACM, 2001.
\newblock Proceedings ISSAC'2001, London, Ontario.

\bibitem{LeLe90}
A.~K. Lenstra and H.~W. {Lenstra, Jr.}
\newblock Algorithms in number theory.
\newblock In J.~van Leeuwen, editor, {\em Handbook of Theoretical Computer
  Science}, volume A: Algorithms and Complexity, chapter~12, pages 674--715.
  North Holland, 1990.

\bibitem{Louboutin02}
St{\'e}phane Louboutin.
\newblock Computation of class numbers of quadratic number fields.
\newblock {\em Math. Comp.}, 71(240):1735--1743 (electronic), 2002.

\bibitem{Mihailescu04}
P.~Mih\u{a}ilescu.
\newblock Fast convolutions meet {M}ontgomery.
\newblock \Preprint, March 2004.

\bibitem{MiAv03}
P.~Mih\u{a}ilescu and R.~Avanzi.
\newblock Efficient quasi-deterministic primality test improving {AKS}.
\newblock Available from {\tt
  http://www-math.uni-paderborn.de/\~{}preda/papers/myaks1.ps}, April 2003.

\bibitem{Morain98a}
F.~Morain.
\newblock Primality proving using elliptic curves: an update.
\newblock In J.~P. Buhler, editor, {\em Algorithmic Number Theory}, volume 1423
  of {\em Lecture Notes in Comput. Sci.}, pages 111--127. Springer-Verlag,
  1998.
\newblock Third International Symposium, ANTS-III, Portland, Oregon, june 1998,
  Proceedings.

\bibitem{Morain04}
F.~Morain.
\newblock La primalit{\'e} en temps polynomial [d'apr\`es {A}dleman, {H}uang~;
  {A}grawal, {K}ayal, {S}axena].
\newblock {\em Ast\'erisque}, pages Exp. No. 917, ix, 205--230, 2004.
\newblock S\'eminaire Bourbaki. Vol. 2002/2003.

\bibitem{Nitaj95}
A.~Nitaj.
\newblock L'algorithme de {C}ornacchia.
\newblock {\em Exposition. Math.}, 13:358--365, 1995.

\bibitem{Schertz02}
R.~Schertz.
\newblock Weber's class invariants revisited.
\newblock {\em J. Th\'eor. Nombres Bordeaux}, 14:325--343, 2002.

\end{thebibliography}
\else
\def\noopsort#1{}\ifx\bibfrench\undefined\def\biling#1#2{#1}\else\def\biling#1%
#2{#2}\fi\def\Inpreparation{\biling{In preparation}{en
  pr{\'e}paration}}\def\Preprint{\biling{Preprint}{pr{\'e}version}}\def\Draft{%
\biling{Draft}{Manuscrit}}\def\Toappear{\biling{To appear}{\`A para\^\i
  tre}}\def\Inpress{\biling{In press}{Sous presse}}\def\Seealso{\biling{See
  also}{Voir {\'e}galement}}\def\Editor{\biling{Ed.}{R{\'e}d.}}

\fi

\end{document}